\documentclass[11pt]{amsart}

\usepackage{amsmath,amssymb,amsthm}
\usepackage{graphicx}
\usepackage{geometry}
\usepackage{bm}
\usepackage{mathtools}
\usepackage{enumitem}
\usepackage{color}
\usepackage{cite}

\newcommand{\R}{\mathbb{R}}
\newcommand{\GammaI}{\Gamma}
\newcommand{\jump}[1]{\left[#1\right]}

\newcommand{\gradS}{\nabla_{\!S}}

\title[A Trace-Based Interface Reduction Method
for Highly Conducting Interfaces]{A Trace-Based Interface Reduction Method \\
for Highly Conducting Interfaces}

\author{So-Hsiang Chou and Do Young Kwak}

\address{Department of Mathematics and Statistics,
Bowling Green State University,
Bowling Green, OH 43403, USA}

\email{chou@bgsu.edu}

\address{Department of Mathematical Sciences,
Korea Advanced Institute of Science and Technology (KAIST),
291 Daehak-ro, Yuseong-gu,
Daejeon 34141, Republic of Korea}

\email{doyoungkwa@kaist.ac.kr}
\date{}

\begin{document}


\begin{abstract}
We develop a reduced interface formulation for elliptic interface problems with highly conducting interfaces. The interface condition consists of continuity of the primal variable together with a jump in the normal flux proportional to the surface Laplacian of the interface trace. Instead of using the solution jump as the interface unknown, we employ the common interface trace and derive a trace-based Schur complement formulation. For prescribed interface trace data, independent extension problems are solved in the two subdomains, leading to a reduced interface equation involving the Dirichlet-to-Neumann jump operator and a surface stiffness operator. Finite-dimensional trace approximations produce compact reduced systems posed only on the interface. Numerical experiments for circular, smooth noncircular, and heart-shaped interfaces illustrate the effectiveness of the method and the role of interface-mode enrichment.
\end{abstract}
\maketitle
\section{Introduction}

Highly conducting thin interphases and imperfect interface models arise
in many applications involving composite materials, diffusion through
thin membranes, and conduction across highly conductive layers.
In such models, the effect of a thin interphase may be represented by
effective interface conditions posed directly on the interface rather
than by explicitly resolving the thin layer itself.
Various imperfect-interface models and interfacial transmission laws have
been studied in mechanics and conduction problems; see, for example,
\cite{MilohBenveniste1999,BenvenisteMilohNeutral,Benveniste2006Interface,GuHe2011}.
In particular, highly conducting interfaces often lead to effective
transmission conditions involving surface differential operators.
Surface-Laplace type interface conditions arise naturally in membrane-type
interface models and in asymptotic descriptions of thin highly conducting
layers.

From the numerical viewpoint, interface problems with surface differential
operators are more challenging than classical transmission problems
because the interface condition itself contains tangential derivatives
along the interface. A variety of fitted, unfitted, immersed, and
nonconforming finite element methods have been developed for elliptic
interface problems; see, for example,
\cite{ChenZou1998,burman2015cutfem,lin2013locking,kwak2010analysis,KwakJin2017}.
Most of these approaches incorporate the interface condition directly
into a global bulk discretization.

The present formulation builds upon our earlier interface-reduction
frameworks for nonlinear and elliptic interface problems; see
\cite{Chou2026FluxRecovery,Chou2026lifting,Chou2026ScalarReduction}.
Those works developed scalar, conservative flux-recovery, and
lifting-based interface reduction methods for classical transmission
problems, where the reduced interface variable is either a scalar
parameter or the solution jump. In contrast, the highly conducting
interface model considered here naturally identifies the common
interface trace as the reduced unknown. This change in perspective
leads to a trace-based Schur complement formulation involving both the
Dirichlet-to-Neumann operator and the surface Laplacian. The resulting
framework extends the interface-reduction philosophy to transmission
problems with surface differential operators and provides a natural
foundation for low-dimensional trace approximations.

The interface conditions considered in this paper take the form
\begin{align}
\jump{u} &= 0\qquad \text{on } \GammaI,
\\
\jump{\beta \partial_n u}
&=
\alpha \Delta_S u_\Gamma\qquad \text{on } \GammaI.
\end{align}
Here $\jump{\cdot}$ denotes the jump across the interface,
$u_\Gamma$ is the common interface trace, and
$\Delta_S$ denotes the surface Laplacian acting tangentially along
$\GammaI$.

These conditions arise as asymptotic limits of thin highly conducting
interphases and appear in the engineering literature on composite media.
Traditional numerical approaches either explicitly resolve thin layers
or solve coupled bulk-interface systems. In contrast, the present work
adopts a reduced interface viewpoint.
The key observation is that since the solution is continuous across the
interface, the natural interface unknown is the common interface trace
rather than the solution jump.

For prescribed interface trace data, independent extension problems are
solved in the two subdomains. The interface condition then induces a
reduced surface equation involving a Dirichlet-to-Neumann jump operator
and a surface Laplacian operator.
The resulting formulation leads naturally to low-dimensional interface
approximations and modal interpretations.

The remainder of the paper is organized as follows.
Section~2 introduces the surface-Laplace interface model and the
associated weak formulation.
Section~3 reviews the tangential differential operators and the surface
integration-by-parts identity used throughout the paper.
Section~4 develops the lifting formulation and derives the reduced
interface equation.
Section~5 introduces the finite-dimensional trace approximation and the
corresponding reduced interface system.
Section~6 discusses spectral representations and geometric trace bases
for the interface approximation space.
Section~7 presents numerical experiments for circular, star-shaped, and
heart-shaped interfaces, including studies of rank enrichment and mesh
refinement.
Finally, Section~8 contains concluding remarks and directions for future
work.

\section{Surface-Laplace Interface Model}

Let $\Omega \subset \R^d$, $d=2,3$, be decomposed into two subdomains
$\Omega^-$ and $\Omega^+$ separated by a smooth interface $\GammaI$.
We consider the elliptic interface problem
\begin{align}\label{eq:model_bulk}
-\nabla \cdot (\beta^{\pm} \nabla u^{\pm}) &= f^{\pm}
\qquad \text{in } \Omega^{\pm},
\\
\jump{u} &= 0
\qquad \text{on } \GammaI,
\\
\jump{\beta \partial_n u}
&=
\alpha \Delta_S u_\Gamma
\qquad \text{on } \GammaI,
\\
\label{eq:model_bc}
u &=0
\qquad \text{on } \partial \Omega.
\end{align}

\subsection{Surface Laplacian}

For a smooth interface $\GammaI$, the surface Laplacian is defined by
\[
\Delta_S v = \gradS \cdot \gradS v.
\]

Here $\nabla_S$ denotes the tangential projection of the ambient gradient
onto the tangent space of the interface.

For a circular interface of radius $R$ parameterized by the angle
$\theta$, one has
\[
\Delta_S v
=
\frac{1}{R^2}
\frac{d^2 v}{d\theta^2}.
\]

Hence Fourier modes diagonalize the surface operator:
\[
\Delta_S e^{ik\theta}
=
-\frac{k^2}{R^2} e^{ik\theta}.
\]

\section{Variational Formulation and Trace-Based Extension Formulation}

Before introducing the reduced trace formulation, it is useful to
observe that the surface-Laplace interface condition naturally leads
to a variational formulation containing both bulk and surface energy
terms.

\subsection{Functional Setting}

The surface-Laplace interface condition introduces a tangential
diffusion term along the interface. Consequently, the interface trace
must possess sufficient regularity for the surface gradient to be
well defined. We therefore consider the energy space
\[
V_\Gamma
=
\Bigl\{
v\in H^1_0(\Omega)
:
v|_\Gamma \in H^1(\Gamma)
\Bigr\}.
\]
Equivalently, writing
$v=(v^-,v^+)$
with
$v^\pm=v|_{\Omega^\pm}$,
we require
\[
v^\pm \in H^1(\Omega^\pm),
\qquad
v^-|_\Gamma=v^+|_\Gamma=:v_\Gamma,
\qquad
v_\Gamma\in H^1(\Gamma).
\]

The weak formulation of
\eqref{eq:model_bulk}
--\eqref{eq:model_bc}
is: find $u\in V_\Gamma$ such that
\begin{equation}
\sum_{\pm}
\int_{\Omega^\pm}
\beta^\pm
\nabla u^\pm\cdot\nabla v^\pm\,dx
+
\alpha
\int_\Gamma
\nabla_S u_\Gamma\cdot\nabla_S v_\Gamma\,ds
=
\sum_{\pm}
\int_{\Omega^\pm}
f^\pm v^\pm\,dx
\label{eq:weak_surface_laplacian}
\end{equation}
for all $v\in V_\Gamma$.

The second term in
\eqref{eq:weak_surface_laplacian}
is a surface stiffness contribution induced by the highly conducting
interface. Thus the surface-Laplace condition introduces tangential
diffusion along the interface.

The associated bilinear form is continuous on $V_\Gamma$.
Under the assumptions
\[
\beta^\pm \ge \beta_0 > 0,
\qquad
\alpha > 0,
\]
standard arguments based on trace inequalities and
Lax--Milgram theory provide a natural framework for
establishing well-posedness of the weak problem.
Since the primary objective of the present paper is the derivation
of a reduced trace formulation, we do not pursue a detailed
functional-analytic analysis here.

This formulation already shows that the problem may be viewed as a
coupled bulk-surface system. The trace reduction method developed
below further reduces the problem to an equation posed entirely on
the interface.
\section{Trace-Based Extension Formulation}

Since the solution is continuous across the interface,
we define the interface trace variable
\[
\gamma := u_\Gamma.
\]

The key idea is to use the interface trace as the reduced variable.
For prescribed trace data $\gamma$, we solve independent extension
problems in the two subdomains.

\subsection{Extension Problems}

Let
\[
H^{\pm} : H^{1/2}(\Gamma)
\to
H^1(\Omega^{\pm})
\]
denote the elliptic extension operators.

For a given interface trace $\gamma\in H^1(\Gamma)$,
the extensions
\[
H^-\gamma,
\qquad
H^+\gamma
\]
are defined as the solutions of
\begin{align}
-\nabla \cdot (\beta^- \nabla (H^-\gamma)) &= 0
\qquad \text{in } \Omega^-,
\\
H^-\gamma &= \gamma
\qquad \text{on } \GammaI,
\\
H^-\gamma &= 0
\qquad \text{on } \partial \Omega \cap \partial \Omega^-,
\end{align}
and
\begin{align}
-\nabla \cdot (\beta^+ \nabla (H^+\gamma)) &= 0
\qquad \text{in } \Omega^+,
\\
H^+\gamma &= \gamma
\qquad \text{on } \GammaI,
\\
H^+\gamma &= 0
\qquad \text{on } \partial \Omega \cap \partial \Omega^+.
\end{align}

Thus the operators $H^{\pm}$ map interface trace data to bulk harmonic
(or more generally elliptic) extensions in the two subdomains.

\subsection{Contribution of the Volume Forcing}

Let
\[
U_f^{\pm}
\]
denote the solutions of the forced bulk problems with homogeneous
interface trace:
\begin{align}
-\nabla \cdot (\beta^- \nabla U_f^-) &= f^-
\qquad \text{in } \Omega^-,
\\
U_f^- &=0
\qquad \text{on } \GammaI,
\\
U_f^- &=0
\qquad \text{on } \partial \Omega \cap \partial \Omega^-,
\end{align}
and
\begin{align}
-\nabla \cdot (\beta^+ \nabla U_f^+) &= f^+
\qquad \text{in } \Omega^+,
\\
U_f^+ &=0
\qquad \text{on } \GammaI,
\\
U_f^+ &=0
\qquad \text{on } \partial \Omega \cap \partial \Omega^+.
\end{align}

The full bulk solutions associated with the trace $\gamma$ are therefore
written as
\[
u^{\pm}(\gamma)
=
H^{\pm}\gamma + U_f^{\pm}.
\]

This decomposition separates the interface contribution from the forcing
contribution and will be useful both analytically and computationally.

\subsection{Dirichlet-to-Neumann Jump Operator}

We now define the Dirichlet-to-Neumann jump operator
\[
S\gamma
:=
\jump{\beta \partial_n(H\gamma)},
\]
where
\[
H\gamma
=
\begin{cases}
H^-\gamma & \text{in } \Omega^-,
\\
H^+\gamma & \text{in } \Omega^+.
\end{cases}
\]

The forcing contribution induces the interface term
\[
g_\Gamma
:=
-\jump{\beta \partial_n U_f}.
\]

Substituting
\[
u^{\pm}(\gamma)
=
H^{\pm}\gamma + U_f^{\pm}
\]
into the interface condition gives
\[
S\gamma - \alpha \Delta_S \gamma
=
g_\Gamma.
\]

The reduced interface equation therefore takes the form
\begin{equation}
S\gamma - \alpha \Delta_S \gamma
=
g_\Gamma.
\label{eq:reduced_interface_equation}
\end{equation}

Equation \eqref{eq:reduced_interface_equation} is the central reduced
interface equation for the highly conducting interface model.

The formulation shows that the original bulk-interface problem can be
reduced to an equation posed entirely on the interface.

\section{Finite-Dimensional Trace Approximation}

Let
\[
\Lambda_m^\gamma
=
\operatorname{span}
\{\psi_1,\dots,\psi_m\}
\subset H^1(\GammaI).
\]
Note that the weak reduced interface equation
\eqref{eq:interface_equation}
contains tangential gradient terms on the interface.
Consequently, the trace variable naturally belongs to \(H^1(\Gamma)\).

We approximate the interface trace by
\[
\gamma_m
=
\sum_{j=1}^m s_j\psi_j.
\]

The reduced equations are obtained by enforcing
\[
\left\langle
\jump{\beta\partial_n H\gamma_m}
-
\alpha\Delta_S\gamma_m,
\psi_i
\right\rangle_{\GammaI}
=
\langle g_\Gamma,\psi_i\rangle_{\GammaI},
\qquad i=1,\dots,m.
\]

Using tangential integration by parts on the interface gives
\begin{equation}\label{eq:interface_equation}
\left\langle
\jump{\beta\partial_n H\gamma_m},
\psi_i
\right\rangle_{\GammaI}
+
\alpha
\left\langle
\gradS\gamma_m,
\gradS\psi_i
\right\rangle_{\GammaI}
=
\langle g_\Gamma,\psi_i\rangle_{\GammaI}.
\end{equation}

Substituting
\[
\gamma_m
=
\sum_{j=1}^m s_j\psi_j
\]
leads to the reduced linear system
\begin{equation}
(S_m+\alpha L_m)s=b_m,
\label{eq:reduced_matrix_system}
\end{equation}
where
\[
(S_m)_{ij}
=
\left\langle
\jump{\beta\partial_n H\psi_j},
\psi_i
\right\rangle_{\GammaI},
\]
and
\[
(L_m)_{ij}
=
\left\langle
\gradS\psi_j,
\gradS\psi_i
\right\rangle_{\GammaI}.
\]

The right-hand side vector is given by
\[
(b_m)_i
=
\langle g_\Gamma,\psi_i\rangle_{\GammaI}.
\]

Equation~\eqref{eq:reduced_matrix_system} is the finite-dimensional
reduced interface system associated with the trace approximation space
$\Lambda_m^\gamma$.

The trace-reduction formulation itself is independent of the particular
choice of interface basis functions.
For general geometries, one may use standard piecewise polynomial finite
element basis functions defined directly on the interface $\GammaI$.
This is the natural choice for complicated interfaces and for
three-dimensional problems.

For circular interfaces, however, Fourier modes provide a particularly
convenient spectral basis because the surface Laplacian diagonalizes in
that representation.
In that case one may choose
\[
\psi_1 = 1,
\qquad
\psi_{2k} = \cos(k\theta),
\qquad
\psi_{2k+1} = \sin(k\theta),
\]
for $k=1,2,\dots$.

If Fourier modes up to frequency $K$ are retained, then the resulting
trace space has dimension
\[
m = 2K+1.
\]
Thus the symbol $m$ denotes the dimension of the reduced trace space,
while $k$ denotes the Fourier frequency index. Increasing $m$ therefore corresponds to enriching the tangential
resolution of the interface approximation.
\section{Modal Interpretation and Interface Enrichment}

The reduced trace formulation admits a natural modal interpretation.
The surface-Laplace operator acts only in the tangential direction
along the interface, and therefore the geometry and parametrization
of the interface play an important role.

\subsection{Circular Interfaces}

For a circular interface of radius $R$, parameterized by the angular
variable $\theta$, the surface Laplacian reduces to
\[
\Delta_S
=
\frac{1}{R^2}
\frac{d^2}{d\theta^2}.
\]

In this case, Fourier modes
\[
\psi_k(\theta)=e^{ik\theta}
\]
diagonalize the surface operator:
\[
\Delta_S \psi_k
=
-\frac{k^2}{R^2}\psi_k.
\]

If the bulk geometry and coefficients are also radially symmetric,
then the Dirichlet-to-Neumann jump operator is diagonal in the same
basis. The reduced equations therefore decouple mode by mode:
\[
\left(
\sigma_k + \alpha \frac{k^2}{R^2}
\right)s_k
=
b_k.
\]
The corresponding eigenvalues of the Dirichlet-to-Neumann jump operator
are denoted by $\sigma_k$.
This interpretation shows that increasing the interface rank corresponds
to enriching tangential interface modes.

\subsection{General Polar Interfaces}

The modal viewpoint is not restricted to circular geometry.
For a general smooth interface represented in polar form
\[
r=r(\theta),
\]
the interface may still be parameterized by the angular variable
$\theta$.

Examples include star-shaped and heart-shaped interfaces such as
\[
r(\theta)=1+0.2\cos(5\theta)
\]
and
\[
r(\theta)=1-\sin(\theta),
\]
which will be used later in the numerical experiments.

In such cases the surface Laplacian is no longer diagonalized exactly
by Fourier modes, since the arclength metric varies along the interface.
Nevertheless, Fourier-type basis functions continue to provide natural
interface approximation spaces.

This observation has important computational consequences.
The geometry enters naturally through the interface parameterization and
the associated tangential derivatives.

For smooth star-shaped interfaces, Fourier enrichment often remains
highly effective because the dominant interface response is still well
captured by low-order tangential modes.

\section{Numerical Experiments}

\subsection{Circular Interface and Modal Rank Enrichment}

We first consider the model problem with a circular interface
\[
\Gamma=
\{(x,y): x^2+y^2=R^2\},
\]
where $R>0$ is fixed.
The circular geometry provides the simplest setting for illustrating the
trace-reduction formulation and the role of tangential interface modes.

The bulk finite element spaces are constructed independently in the two
subdomains, allowing the interface coupling conditions to be imposed
through the reduced interface formulation.
The emphasis of the present work is not the construction of a new bulk
interface discretization, but rather the low-dimensional interface reduction
associated with the interface unknowns.

For the circular geometry, the interface admits a natural Fourier
representation in the angular variable $\theta$.
Consequently, the reduced interface space is enriched progressively by
trigonometric modes, allowing us to examine how the approximation error
depends on the number of retained interface modes.
For a circle parameterized by the angular variable $\theta$, the surface
Laplacian reduces to
\[
\Delta_S
=
\frac{1}{R^2}
\frac{d^2}{d\theta^2}.
\]
Consequently, Fourier modes diagonalize the surface operator:
\[
\Delta_S e^{ik\theta}
=
-
\frac{k^2}{R^2}
 e^{ik\theta}.
\]

This structure makes the circular geometry particularly useful for
understanding the interaction between the bulk extension operator and the
surface diffusion operator.

In the numerical experiments below, the interface trace is approximated by
real Fourier basis functions:
\[
\Lambda_m^\gamma
=
\operatorname{span}
\{1,
\cos\theta,
\sin\theta,
\dots,
\cos(k\theta),
\sin(k\theta)
\}.
\]

The reduced trace approximation takes the form
\[
\gamma_m(\theta)
=
\sum_{j=1}^m s_j \psi_j(\theta).
\]

For each basis function $\psi_j$, independent extension problems are
solved in the interior and exterior subdomains:
\[
H^-\psi_j,
\qquad
H^+\psi_j.
\]
The corresponding flux jump response
\[
[\beta \partial_n H\psi_j]
\]
is then used to assemble the reduced Dirichlet-to-Neumann matrix
\[
(S_m)_{ij}
=
\langle
[\beta \partial_n H\psi_j],
\psi_i
\rangle_\Gamma.
\]

The surface diffusion contribution is assembled independently through
\[
(L_m)_{ij}
=
\int_\Gamma
\nabla_S \psi_j \cdot \nabla_S \psi_i\,ds.
\]

For a circular interface,
\[
\nabla_S
=
\frac1R \frac{d}{d\theta},
\qquad
ds = R\,d\theta,
\]
and therefore
\[
(L_m)_{ij}
=
\frac1R
\int_0^{2\pi}
\psi_j'(\theta)
\psi_i'(\theta)
\,d\theta.
\]

Thus the surface stiffness matrix is assembled entirely from one-dimensional
quadrature along the interface.
No surface finite element mesh or differential-geometry machinery is
required.

The reduced interface system is
\[
(S_m + \alpha L_m)s=b_m.
\]
The important feature of this formulation is that the bulk extension solves
remain completely decoupled from the surface diffusion operator.
The extension problems are standard elliptic solves with prescribed trace
data, while the surface-Laplace contribution enters only through the
low-dimensional interface matrix $L_m$.

To illustrate the modal structure of the reduced system, we consider a
manufactured forcing contribution whose interface flux jump contains only
low-order Fourier modes:
\[
[\beta \partial_n U_f]
=
0.5\cos\theta
+
0.25\sin(2\theta).
\]

Consequently, the exact reduced trace contains only the modes
\[
\cos\theta,
\qquad
\sin(2\theta).
\]

The numerical experiment therefore provides a direct test of the
interface-rank enrichment philosophy.
As the interface rank $m$ increases, the approximation space eventually
contains all active tangential modes present in the forcing.
Once these modes are included, the reduced solution becomes exact up to
machine precision.

Table~\ref{table:circular_rank_enrichment} reports the interface error for
several interface ranks.
A high-rank solution with $m=25$ is used as the reference solution.

\begin{table}[htbp]
\centering
\begin{tabular}{c|cc}
\hline
$m$
& $L^2(\Gamma)$ error
& $L^{\infty}(\Gamma)$ error
\\
\hline
1
& $3.422291\times 10^{-2}$
& $3.015191\times 10^{-2}$
\\
3
& $1.098114\times 10^{-2}$
& $8.761682\times 10^{-3}$
\\
5
& $8.908318\times 10^{-18}$
& $1.561251\times 10^{-17}$
\\
7
& $8.420471\times 10^{-18}$
& $1.387779\times 10^{-17}$
\\
9
& $8.072781\times 10^{-18}$
& $1.214306\times 10^{-17}$
\\
\hline
\end{tabular}
\caption{Rank enrichment for the circular interface problem.}
\label{table:circular_rank_enrichment}
\end{table}

Several important observations follow from Table
\ref{table:circular_rank_enrichment}.

When $m=1$, the trace space contains only the constant mode and therefore
cannot represent the active forcing modes.
The error is correspondingly large.

For $m=3$, the approximation space contains the modes
\[
1,
\cos\theta,
\sin\theta,
\]
so the dominant $\cos\theta$ component is captured.
However, the $\sin(2\theta)$ mode is still missing.
The error therefore decreases, but remains non-negligible.

When $m=5$, the approximation space additionally contains
\[
\cos(2\theta),
\qquad
\sin(2\theta),
\]
which completes the set of active modes appearing in the manufactured
forcing.
At this point the reduced solution becomes exact up to machine precision.
Further increases in the interface rank produce essentially no additional
improvement.

This experiment clearly demonstrates the central principle of the
trace-reduction formulation:
interface enrichment resolves the tangential interface physics while the
bulk extension mechanism remains unchanged.
The geometric and surface-diffusion effects are confined entirely to a
small interface system.

\begin{figure}[htbp]
\centering
\includegraphics[width=0.62\textwidth]{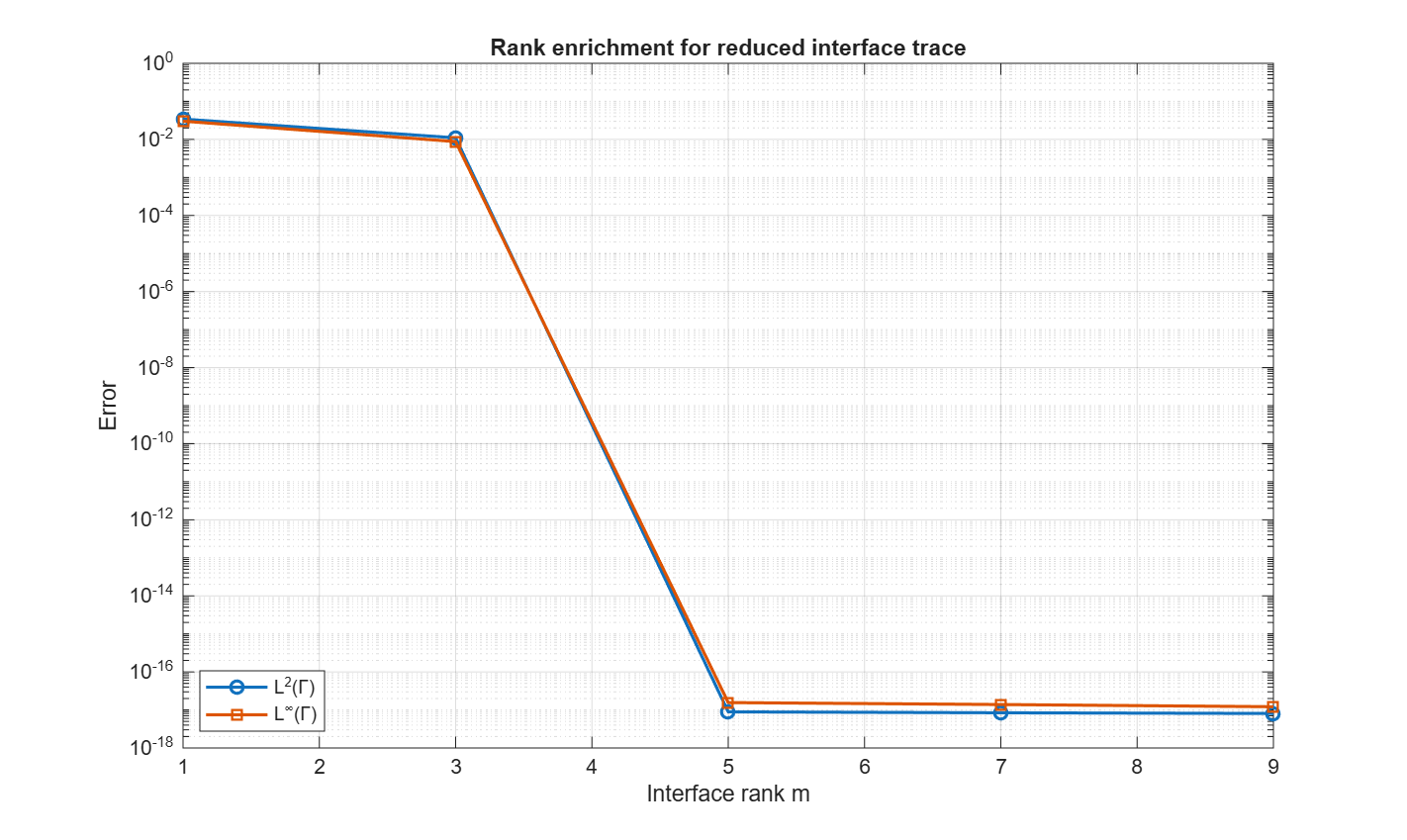}
\caption{Interface-rank enrichment for the circular interface problem.
The error drops to machine precision once all active Fourier modes are included.}
\label{fig:circular_trace_solution}
\end{figure}

Figure~\ref{fig:circular_trace_solution} shows the reduced interface trace
for the case $m=9$.

The circular example provides a transparent illustration of the modal
structure of the surface-Laplace interface condition.
More complicated geometries will be considered in the following sections,
where exact Fourier diagonalization is no longer available.
Nevertheless, the same trace-reduction philosophy continues to apply:
independent bulk extension solves are combined with low-dimensional
interface operators defined entirely along the interface.

\subsection{Star-Shaped Interface}

We next consider a smooth noncircular interface.  This example is meant to
show that the trace-reduction formulation is not tied to the special
Fourier diagonalization available for a circular interface.

Let the interface be given in polar form by
\begin{equation}
\Gamma=
\{(r(\theta)\cos\theta,r(\theta)\sin\theta):0\leq \theta<2\pi\},
\label{eq:star_interface_param}
\end{equation}
where
\begin{equation}
r(\theta)=r_0+\varepsilon \sin(k\theta).
\label{eq:star_radius}
\end{equation}
In the computations below we use
\[
r_0=0.5,
\qquad
\varepsilon=0.08,
\qquad
k=5.
\]
Thus the interface has a five-fold star-shaped perturbation of a circle.
The interface remains smooth and strictly star-shaped with respect to the
origin.

The main difference from the circular case is that arclength is no longer
proportional to the angular variable.  Instead,
\begin{equation}
ds
=
q(\theta)\,d\theta,
\qquad
q(\theta)
=
\sqrt{r(\theta)^2+r'(\theta)^2}.
\label{eq:polar_arclength_metric}
\end{equation}
Consequently, tangential differentiation is given by
\begin{equation}
\frac{d}{ds}
=
\frac{1}{q(\theta)}
\frac{d}{d\theta}.
\label{eq:tangential_derivative_polar}
\end{equation}

Therefore, for interface basis functions $\psi_i(\theta)$ and
$\psi_j(\theta)$, the surface stiffness matrix becomes
\begin{equation}
(L_m)_{ij}
=
\int_\Gamma
\nabla_S\psi_j\cdot \nabla_S\psi_i\,ds
=
\int_0^{2\pi}
\frac{\psi_j'(\theta)\psi_i'(\theta)}
{q(\theta)}\,d\theta.
\label{eq:polar_surface_stiffness}
\end{equation}
This formula is the only change needed in the surface part of the
reduced system.  The bulk extension problems are unchanged in principle:
for each trace basis function $\psi_j$, we solve the independent extension
problems in the two subdomains and compute the corresponding flux jump
response.

This example therefore illustrates one of the main advantages of the
trace-extension viewpoint.  The geometric complexity of the curved
interface is confined to one-dimensional quadrature along $\Gamma$ and to
the construction of the bulk extension problems.  The surface-Laplace term
is not assembled as part of a global bulk-surface finite element system.
It enters only through the reduced matrix $L_m$.

As in the circular example, we use real Fourier functions of the parameter
$\theta$ as the interface approximation basis:
\[
1,
\cos\theta,
\sin\theta,
\cos(2\theta),
\sin(2\theta),
\dots.
\]
These functions no longer diagonalize the surface Laplacian, since the
metric factor $q(\theta)$ is not constant.  Nevertheless, they provide a
natural global approximation space for smooth star-shaped interfaces.

For each interface rank $m$, we assemble
\begin{equation}
(S_m+\alpha L_m)s=b_m.
\label{eq:star_reduced_system}
\end{equation}
Here $S_m$ is obtained from the bulk extension responses and $L_m$ is
computed from \eqref{eq:polar_surface_stiffness}.  The right-hand side
contains the contribution from the manufactured forcing part $U_f$.

The numerical experiment is designed to test whether interface enrichment
continues to be effective when exact Fourier diagonalization is lost.  A
high-rank reduced solution is used as the reference solution.  The error is
then measured for a sequence of interface spaces of increasing dimension.

\begin{table}[htbp]
\centering
\begin{tabular}{c|cc}
\hline
$m$
& $L^2(\Gamma)$ error
& $L^{\infty}(\Gamma)$ error
\\
\hline
1  & $1.230135\times 10^{-1}$ & $1.014424\times 10^{-1}$ \\
3  & $1.228736\times 10^{-1}$ & $9.968503\times 10^{-2}$ \\
5  & $1.228733\times 10^{-1}$ & $9.969220\times 10^{-2}$ \\
7  & $1.228661\times 10^{-1}$ & $9.961740\times 10^{-2}$ \\
9  & $1.227320\times 10^{-1}$ & $1.001117\times 10^{-1}$ \\
13 & $2.470175\times 10^{-2}$ & $2.638302\times 10^{-2}$ \\
17 & $2.475577\times 10^{-2}$ & $2.634641\times 10^{-2}$ \\
\hline
\end{tabular}
\caption{Rank enrichment for the star-shaped interface problem.}
\label{table:star_rank_enrichment}
\end{table}

The expected behavior is different from the circular case.  In the circular
experiment, the manufactured data contained only finitely many Fourier
modes and the reduced solution became exact once those modes were included.
For the star-shaped interface, the variable metric factor $q(\theta)$
couples Fourier modes through the surface stiffness matrix.  As a result,
one should expect a more gradual decay of the error as $m$ increases.

This slower decay is not a weakness of the method.  Rather, it reflects the
fact that the geometry itself transfers information among tangential modes.
The additional enrichment remains localized to the interface.
The same bulk extension machinery is reused, while the
additional geometric and surface-diffusion effects are represented in the
small reduced system \eqref{eq:star_reduced_system}.

\begin{figure}[htbp]
\centering
\includegraphics[width=0.62\textwidth]{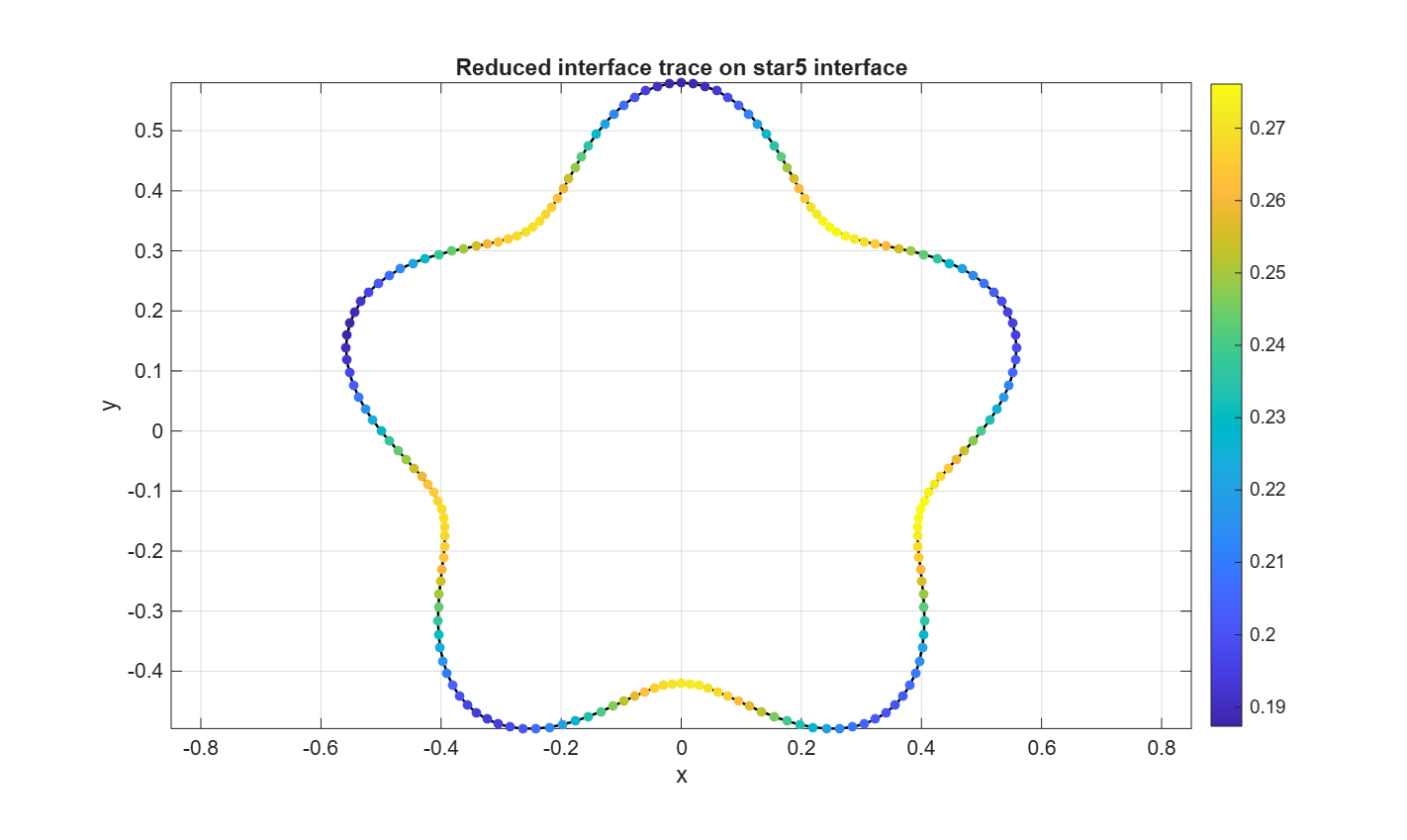}
\caption{
Computed reduced interface trace on the five-petal star-shaped interface.
The color represents the reduced trace value $\gamma_m$ along $\Gamma$.
}
\label{fig:star_trace_solution}
\end{figure}

Figure~\ref{fig:star_trace_solution} illustrates the reduced interface trace directly on the
star-shaped geometry.  The trace variable is represented entirely along
the interface, while the bulk extension problems remain decoupled.

This example provides the bridge between the analytically transparent
circular case and the more geometrically challenging interfaces considered next. It demonstrates that the trace-reduction formulation continues to
be straightforward even when the surface Laplacian is no longer diagonal
in the chosen basis.

\subsection{Heart-Shaped Interface}

We now consider a more challenging heart-shaped interface.
The purpose of this example is to examine whether the trace-reduction
philosophy remains effective for interfaces with more complicated
geometry than the smooth star-shaped interface considered previously.

The interface is given parametrically by
\begin{equation}
\Gamma=
{(x(\theta),y(\theta)):0\leq \theta<2\pi},
\label{eq:heart_interface_param}
\end{equation}
where
\begin{equation}
x(\theta)=r(\theta)\cos\theta,
\qquad
y(\theta)=r(\theta)\sin\theta,
\label{eq:heart_parametric_xy}
\end{equation}
and
\begin{equation}
r(\theta)=1-\sin\theta.
\label{eq:heart_radius}
\end{equation}

As before, the arclength metric is
\begin{equation}
q(\theta)=\sqrt{r(\theta)^2+r'(\theta)^2},
\qquad
\,ds=q(\theta)\,d\theta,
\label{eq:heart_metric}
\end{equation}
and tangential differentiation satisfies
\begin{equation}
\frac{d}{ds}=
\frac1{q(\theta)}
\frac{d}{d\theta}.
\label{eq:heart_tangential_derivative}
\end{equation}
Consequently, the surface stiffness matrix is assembled through
\begin{equation}
(L_m)_{ij}=\int_0^{2\pi}
\frac{\psi_j'(\theta)\psi_i'(\theta)}{q(\theta)}
\,d\theta.
\label{eq:heart_surface_matrix}
\end{equation}

Compared with the star-shaped interface, the present geometry produces
stronger localized variation in the arclength metric $q(\theta)$.
Nevertheless, the bulk extension problems are again solved independently
in the two subdomains. Thus the geometric complexity enters only through
the interface quadrature and the extension geometry, while the
surface-Laplace contribution remains confined to the reduced interface
matrix.

For the interface approximation space, we again employ global
Fourier-type basis functions in the parameter variable $\theta$:
\[
1,\quad
\cos\theta,\quad
\sin\theta,\quad
\cos(2\theta),\quad
\sin(2\theta),\dots .
\]
Although exact modal diagonalization is no longer available, these basis
functions continue to provide an effective global representation of the
interface trace.

The reduced system has the same algebraic form as before:
\begin{equation}
(S_m+\alpha L_m)s=b_m.
\label{eq:heart_reduced_system}
\end{equation}

Since tangential differentiation is scaled by
\[
\nabla_S=\frac1{q(\theta)}\frac{d}{d\theta},
\]
the geometry enters the reduced surface operator through the metric
factor \(q(\theta)\).

Despite this additional complexity, the numerical results below show
that interface-rank enrichment continues to reduce the error effectively.
This example illustrates that increased geometric complexity does not
require a fundamentally different bulk discretization strategy. Traditional interface methods would
typically address such a geometry through additional local mesh
refinement. In contrast, the trace-reduction formulation seeks to
capture the additional complexity primarily through interface enrichment.

A high-rank reduced solution is again used as the reference solution.
The reduced trace is then computed for a sequence of interface ranks.

\begin{table}[htbp]
\centering
\begin{tabular}{c|cc}
\hline
$m$
& $L^2(\Gamma)$ error
& $L^{\infty}(\Gamma)$ error
\\
\hline
1  & $1.163554\times 10^{-2}$ & $1.498769\times 10^{-2}$ \\
3  & $9.174113\times 10^{-3}$ & $1.164148\times 10^{-2}$ \\
5  & $4.398064\times 10^{-2}$ & $4.594850\times 10^{-2}$ \\
7  & $6.303609\times 10^{-3}$ & $7.720991\times 10^{-3}$ \\
9  & $5.449365\times 10^{-3}$ & $7.851867\times 10^{-3}$ \\
13 & $6.541914\times 10^{-4}$ & $7.665676\times 10^{-4}$ \\
17 & $3.065561\times 10^{-4}$ & $4.370653\times 10^{-4}$ \\
21 & $1.594905\times 10^{-4}$ & $2.417035\times 10^{-4}$ \\
25 & $1.012410\times 10^{-4}$ & $1.752753\times 10^{-4}$ \\
31 & $1.002378\times 10^{-4}$ & $1.422356\times 10^{-4}$ \\
37 & $7.812372\times 10^{-5}$ & $7.780763\times 10^{-5}$ \\
\hline
\end{tabular}
\caption{Rank enrichment for the heart-shaped interface problem.}
\label{table:heart_rank_enrichment}
\end{table}

The results indicate that the reduced trace remains well approximated by
interface enrichment even for this more demanding geometry. Although the
convergence is less structured than in the circular case, the error
continues to decrease as the interface rank increases.

\begin{figure}[htbp]
\centering
\includegraphics[width=0.62\textwidth]{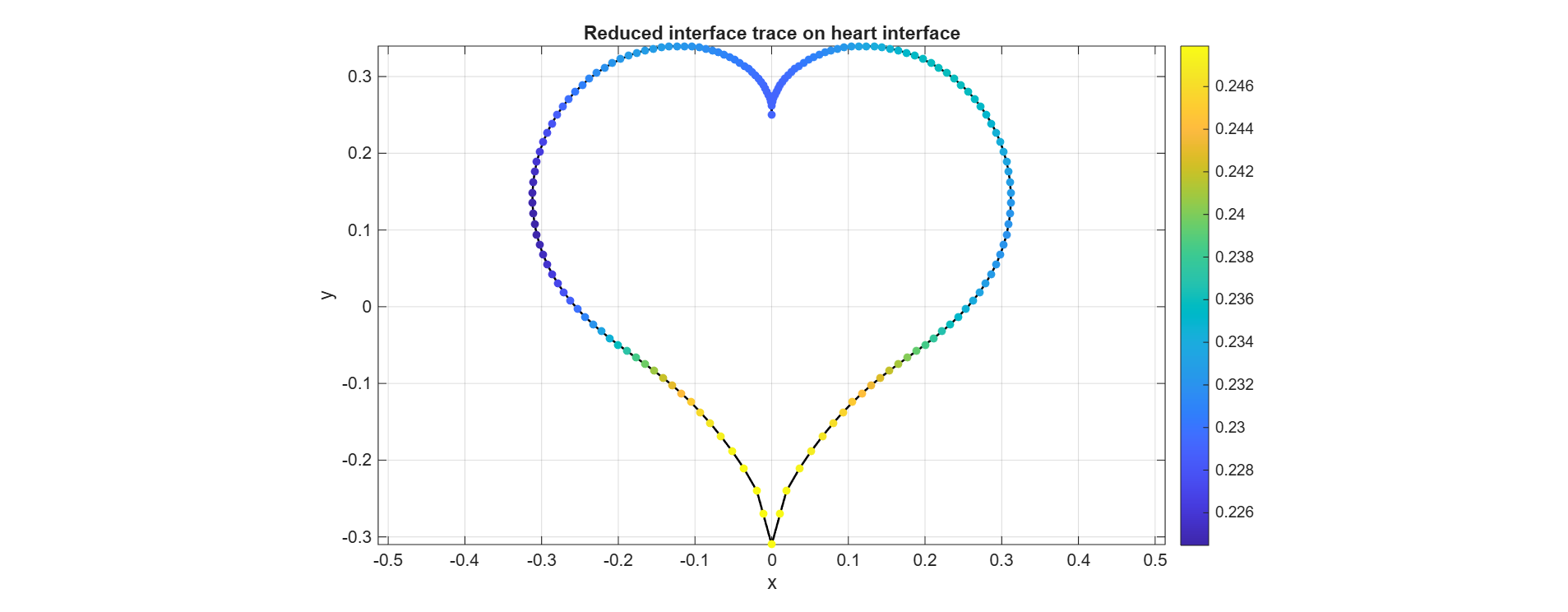}
\caption{
Computed reduced interface trace on the heart-shaped interface.
The color represents the reduced trace value $\gamma_m$ along $\Gamma$.
}
\label{fig:heart_trace_solution}
\end{figure}

Figure~\ref{fig:heart_trace_solution} illustrates the reduced interface
trace directly on the heart-shaped interface. The trace variable is
represented entirely along the interface, while the bulk extension
problems remain decoupled.

The heart-shaped example illustrates that the trace-reduction framework
remains computationally manageable even for geometrically challenging
interfaces. The complexity of the geometry is absorbed primarily into the
low-dimensional interface operators rather than into repeated global bulk
refinement. The bulk extension problems themselves remain standard
elliptic solves.

\subsection{Rank Enrichment versus Mesh Refinement}

One of the main motivations behind the present formulation is the
observation that, for interface problems with surface diffusion effects,
the dominant numerical difficulty often resides on the interface itself
rather than in the bulk discretization.
The purpose of this subsection is therefore to compare two different
strategies for improving accuracy:

\begin{itemize}
\item bulk mesh refinement with fixed interface rank;
\item interface-rank enrichment with a fixed bulk mesh.
\end{itemize}

The comparison is intended to clarify where the principal approximation
error originates in the trace-reduction formulation.

Throughout this experiment we use the star-shaped interface introduced in
Section~7.2.
The bulk meshes are generated independently in the interior and exterior
subdomains, while the reduced interface approximation uses the Fourier
trace basis
\[
1,
\cos\theta,
\sin\theta,
\cos(2\theta),
\sin(2\theta),
\dots.
\]

The reduced system has the form
\begin{equation}
(S_m+\alpha L_m)s=b_m,
\label{eq:rank_mesh_compare_system}
\end{equation}
where $S_m$ is obtained from the bulk extension responses and $L_m$
represents the surface-Laplace contribution.

\subsubsection*{Fixed-rank mesh refinement}

We first keep the interface rank fixed and refine only the bulk meshes.
Specifically, we use a fixed interface rank
\[
m=5,
\]
while decreasing the mesh size parameter $H_{\max}$.

A high-rank computation on a very fine mesh is used as the reference
solution.
Table~\ref{table:fixed_rank_mesh_refine} reports the resulting interface
errors.

\begin{table}[htbp]
\centering
\begin{tabular}{c|cc}
\hline
$H_{\max}$
& $L^2(\Gamma)$ error
& $L^\infty(\Gamma)$ error
\\
\hline
0.090 & $9.386265\times 10^{-2}$ & $5.879512\times 10^{-2}$ \\
0.070 & $9.333609\times 10^{-2}$ & $5.850061\times 10^{-2}$ \\
0.055 & $9.163721\times 10^{-2}$ & $5.776560\times 10^{-2}$ \\
0.040 & $9.278833\times 10^{-2}$ & $5.816984\times 10^{-2}$ \\
\hline
\end{tabular}
\caption{Fixed-rank mesh refinement with $m=5$.}
\label{table:fixed_rank_mesh_refine}
\end{table}

The error decreases only mildly under bulk refinement.
The results suggest that the dominant approximation error is no longer
primarily caused by insufficient bulk resolution.
Instead, the principal limitation arises from the restricted interface
approximation space.

In other words, once the extension problems are adequately resolved,
further bulk refinement alone cannot efficiently recover missing
tangential interface modes.

\subsubsection*{Fixed-mesh rank enrichment}

We next reverse the strategy.
The bulk mesh is now kept fixed, while the interface rank $m$ is
increased.

Table~\ref{table:fixed_mesh_rank_enrich} reports the corresponding
errors.

\begin{table}[htbp]
\centering
\begin{tabular}{c|cc}
\hline
$m$
& $L^2(\Gamma)$ error
& $L^\infty(\Gamma)$ error
\\
\hline
1  & $1.688396\times 10^{-1}$ & $1.014495\times 10^{-1}$ \\
3  & $1.686639\times 10^{-1}$ & $9.973528\times 10^{-2}$ \\
5  & $1.686636\times 10^{-1}$ & $9.974311\times 10^{-2}$ \\
7  & $1.686539\times 10^{-1}$ & $9.967854\times 10^{-2}$ \\
9  & $1.684947\times 10^{-1}$ & $1.001468\times 10^{-1}$ \\
13 & $3.503522\times 10^{-2}$ & $2.639523\times 10^{-2}$ \\
17 & $3.510978\times 10^{-2}$ & $2.634788\times 10^{-2}$ \\
\hline
\end{tabular}
\caption{Fixed-mesh interface-rank enrichment.}
\label{table:fixed_mesh_rank_enrich}
\end{table}

In contrast with the previous experiment, the error now decreases much
more rapidly as the interface rank increases.
This behavior demonstrates that the remaining unresolved structure is
primarily localized on the interface.

The numerical results therefore support the main computational principle
of the present work:
after the bulk extension operator has been adequately approximated,
additional accuracy is more effectively obtained through interface
enrichment rather than repeated global mesh refinement.

\begin{figure}[htbp]
\centering

\begin{minipage}{0.48\textwidth}
\centering
\includegraphics[width=\textwidth]
{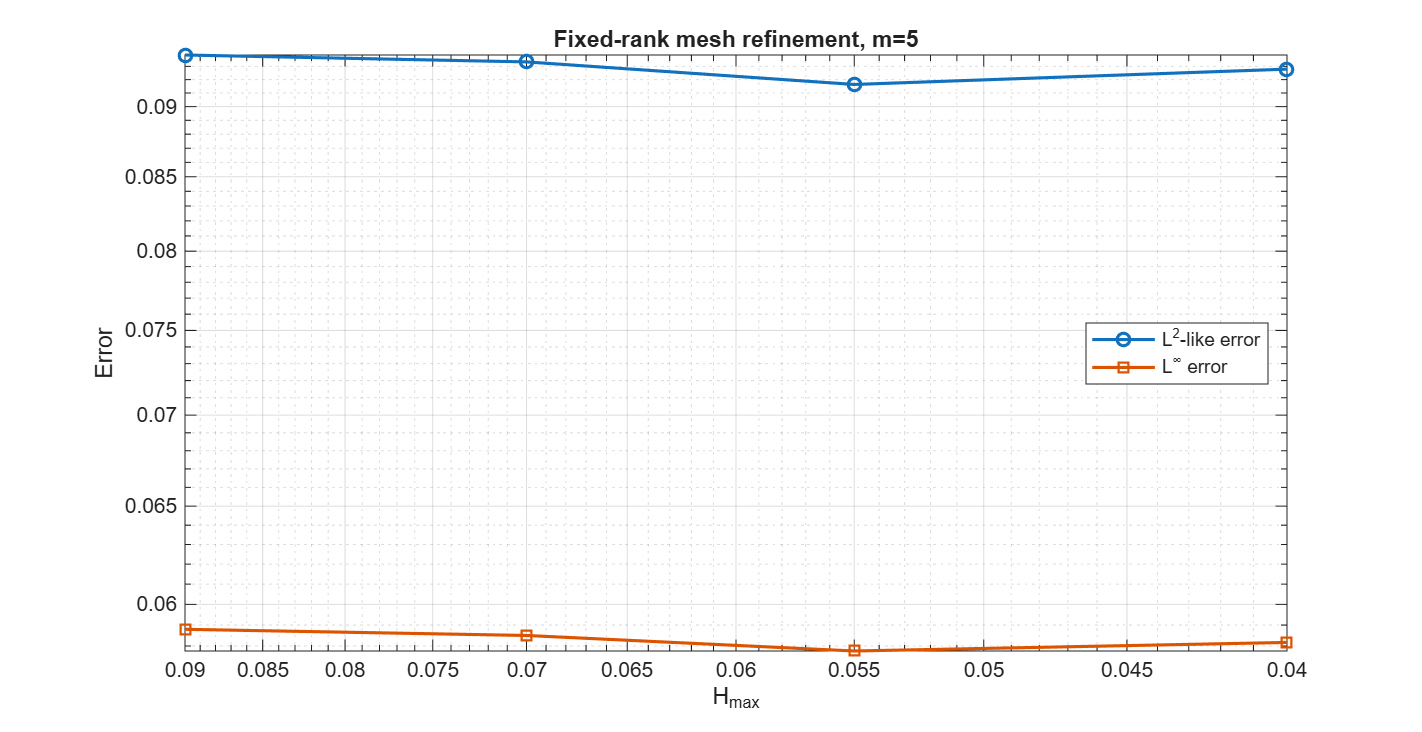}
\end{minipage}
\hfill
\begin{minipage}{0.48\textwidth}
\centering
\includegraphics[width=\textwidth]
{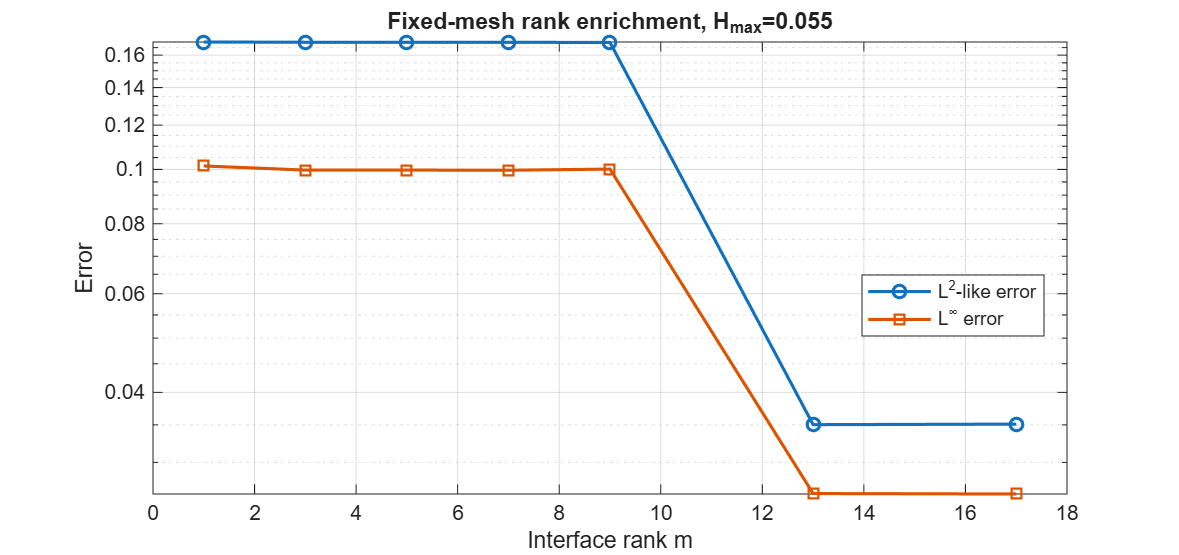}
\end{minipage}

\caption{
Comparison between bulk mesh refinement and interface-rank enrichment.
Left: fixed-rank mesh refinement with $m=5$ produces only marginal
improvement. Right: fixed-mesh interface enrichment produces substantial
error reduction once the active interface modes are resolved.
}
\label{fig:rank_vs_mesh}
\end{figure}

\subsubsection*{Error reduction mechanisms}

The previous experiments highlight an important distinction between
traditional bulk finite element refinement and the present
trace-reduction formulation. Classical finite element methods typically improve accuracy through
global refinement of the bulk approximation space.
In the present formulation, however, the solution is decomposed into
\[
u^\pm(\gamma)
=
H^\pm\gamma
+
U_f^\pm.
\]
Here the bulk extension operators $H^\pm$ are reused throughout the
computation, while the remaining approximation error is largely governed
by the interface trace variable $\gamma$. Consequently, the dominant error mechanism is often associated with the
resolution of tangential interface structure rather than with the bulk
elliptic solves themselves. From this viewpoint, interface-rank enrichment plays a role analogous to
a localized spectral enrichment on the interface.
The complexity of the geometry and the surface-Laplace operator is
absorbed into the low-dimensional interface system
\eqref{eq:rank_mesh_compare_system},
while the bulk solves remain standard elliptic subproblems.

This perspective contrasts with the usual strategy of uniform
bulk mesh refinement.
Instead of repeatedly refining the entire computational domain, the
trace-reduction formulation attempts to identify and enrich the lower-
dimensional structure where the dominant unresolved physics resides.

The principal computational advantage of the present
formulation is that, after the bulk extension responses
have been computed, the interface unknown is determined
from a reduced system of dimension $m$.
For the numerical experiments reported here,
$m$ ranges from 5 to 37, whereas the corresponding
bulk finite element discretizations contain several
thousand degrees of freedom.
Thus the surface-Laplace coupling is represented by
a low-dimensional interface problem whose size is
independent of the bulk mesh resolution.
\section{Concluding Remarks}

We introduced a trace-based interface reduction formulation for
surface-Laplace interface conditions associated with highly conducting
interfaces. The formulation uses the common interface trace as the reduced unknown and combines independent bulk extension problems with a reduced surface interface equation. The resulting reduced operator consists of a Dirichlet-to-Neumann jump map together with a surface stiffness operator.
The numerical experiments indicate that, once the bulk extension
operators are adequately resolved, the dominant approximation error is
often governed by the interface trace space rather than by the bulk mesh
itself.

Future work includes extensions to coupled interface systems,
Stokes and elasticity problems with surface operators,
and eigenvalue problems involving interface dynamics.

\bibliographystyle{abbrv}
\bibliography{references}

\end{document}